\documentclass[graybox]{svmult}

\usepackage{type1cm}        
\usepackage{makeidx}         
\usepackage{graphicx}        
\usepackage{multicol}        
\usepackage[bottom]{footmisc}

\usepackage[T1]{fontenc}
\usepackage{newtxtext}       %
\usepackage[varvw]{newtxmath}       

\usepackage{hhline}
\usepackage{comment}

\makeindex             

\def\PP{{{\rm l}\kern - .15em {\rm P} }}
\def\PN2{{\PP_{N}-\PP_{N-2}}}

\newcommand{\cO}{\mathcal{O}}

\newcommand{\bphi}{\boldsymbol{\varphi}}

\newcommand{\bu}{\boldsymbol{u}}

\newcommand{\bur}{{\boldsymbol{u}}_r}

\newcommand{\apr}{{\it a priori} error bounds}
\newcommand{\apo}{{\it a posteriori} error bounds}

\begin{document}

\title*{{\it A Priori} Error Bounds for POD-ROMs for Fluids: A Brief Survey}
\author{Francesco Ballarin\orcidID{0000-0001-6460-3538} and\\ Traian Iliescu\orcidID{0000-0003-1437-7362}}
\institute{
Francesco Ballarin \at Università Cattolica del Sacro Cuore, Department of Mathematics and Physics, 25133 Brescia, Italy, \email{francesco.ballarin@unicatt.it}
\and Traian Iliescu \at Virginia Tech, Department of Mathematics, Blacksburg, VA 24061, USA, \email{iliescu@vt.edu}}
%
%
\maketitle

\abstract{
Galerkin reduced order models (ROMs), e.g., based on proper orthogonal decomposition (POD) or reduced basis methods, have achieved significant success in the numerical simulation of fluid flows.
The ROM numerical analysis, however, is still being developed. 
In this paper, we take a step in this direction and present a survey of \apr, 
with a particular focus on POD-based ROMs.
Specifically, we outline the main components of ROM \apr,  emphasize their practical importance, and discuss significant contributions to \apr \ for ROMs for fluids. 
}

\section{Introduction}
\label{sec:introduction}


Galerkin reduced order models (ROMs, in the following) leverage available data from full order models (FOMs, in the following) to construct a data-driven basis. 
While in principle any FOM can be used, in this review will we focus on FOMs based on the finite element method (FEM), which are then used to provide data for the basis construction.
Such a basis is typically built by means of the proper orthogonal decomposition (POD) algorithm or the reduced basis method (RBM) during the so-called offline phase.
The constructed basis is then used within a Galerkin method to query a model with significantly reduced dimensionality during the so-called online phase. The dimensionality of the ROM is typically much lower than 
the dimensionality of the corresponding FOM.
ROMs have been used successfully in providing an efficient and reasonably accurate solution for the numerical simulation of fluid flows, particularly in the laminar regime.


The numerical analysis of ROMs for fluids is divided into two primary branches, depending on the method employed during the offline phase:
(i) numerical analysis for POD, and
(ii) numerical analysis for RBM.
These branches, though they began around the same time, have evolved in {\it fundamentally different directions}: Indeed, the numerical analysis of RBM has concentrated on proving \apo, which are instrumental in constructing suitable bases for parameter-dependent problems.
On the other hand, the numerical analysis of POD has primarily focused on proving \apr.
It is worth noting that the numerical analysis of RBM is relatively well-developed, both for general problems and for fluids~\cite{haasdonk2017reduced,hesthaven2015certified,quarteroni2015reduced}.
However, the numerical analysis of POD for fluids is still being developed. 

In this paper, we present a brief survey of the numerical analysis of POD-based ROMs for fluids flows, focusing on \apr.
In Section~\ref{sec:a-priori-posteriori-error-bounds}, we outline the main components of \apr \ and 
include a brief qualitative comparison with \apo.
In Section~\ref{sec:practical-importance}, we emphasize the  importance of the \apr \ in practical settings.
In Section~\ref{sec:major-contributions}, we discuss major contributions made to the development of \apr \ for ROMs of fluid flows, focusing on the most recent results.
Finally, in Section~\ref{sec:conclusions}, we draw conclusions and outline future research directions.

\section{{\it A Priori} vs. {\it A Posteriori} Error Bounds}
    \label{sec:a-priori-posteriori-error-bounds}

In this review paper, we focus on the incompressible Navier-Stokes equations (NSE) on a bounded domain $\Omega \subset \mathbb{R}^d$, with $d = 2, 3$. Let $\bu$ denote the exact velocity of the NSE, $\bu_h$ the FEM velocity approximation, and $\bur$ the ROM velocity approximation. We will not introduce a symbol for the pressure unknown since, as typically done in ROMs built with a POD basis, we assume that the basis is already weakly divergence-free, and hence the pressure drops out from the ROM system. We refer the reader to alternative approaches described, e.g., in \cite{ballarin2015supremizer,rozza2007stability,stabile2018finite} and the references discussed in Section~\ref{sec:pressure} if online pressure recovery is desired.

\subsection{{\it A Priori} Error Bounds}
    \label{sec:a-priori-error-bounds}


The ROM \apr \ are generally of the following form:
\begin{eqnarray}
    \| \bu - \bur \|
    \leq C \ f \left(\Delta t, h, r, \{\lambda_i\}_{i=r+1}^{R}, \{\bphi_i\}_{i=r+1}^{R} \right), 
    \label{eqn:a-priori-error-bound}
\end{eqnarray}
where:
\begin{itemize}
\item $\| \cdot \|$ denotes, e.g., the $L^2(\Omega)$ norm,
\item $\Delta t$ is the time step size of the FOM and ROM time discretizations (although different time steps can be considered for the FOM and ROM time discretizations), 
\item $h$ is the mesh size of the FOM spatial discretization,
\item $r$ is the dimension of the ROM basis,
\item $R$ is the rank of the snapshot matrix used to construct the ROM basis,
\item $\{\lambda_i\}_{i=r+1}^{R}$ are eigenvalues in the POD eigenvalue problem,  
and $\{\bphi_i\}_{i=r+1}^{R}$ are POD basis functions,
\item $f$ is a function that describes the \apr, its arguments representing the FOM and ROM model parameters listed above,
\item $C$ is a generic constant that can depend on the problem data (e.g., domain, initial conditions, boundary conditions, and exact solution, $\bu$), but not on the FOM and ROM parameters.
\end{itemize}

\noindent We remark the following about the {\it a priori} error bound~\eqref{eqn:a-priori-error-bound}:
\begin{enumerate}
\item The {\it a priori} error bound~\eqref{eqn:a-priori-error-bound} is {\it asymptotic}:
It tells us how the error bound behaves in the asymptotic regime of the FOM and ROM parameters, i.e., for
$h \rightarrow 0$, $\Delta t \rightarrow 0$, and $r \rightarrow R$.
\item 
The {\it a priori} error bound~\eqref{eqn:a-priori-error-bound} generally yields the ROM's {\it rate of convergence}, i.e., it tells us how fast the error goes to $0$ as $h \rightarrow 0$, $\Delta t \rightarrow 0$, and $r \rightarrow R$.
Specifically, the {\it a priori} error bound~\eqref{eqn:a-priori-error-bound} yields rates of convergence that are generally expressed by using the celebrated {\it ``big O'' notation}, which is central in the numerical analysis of classical numerical methods, e.g., 
the FEM.
For example, the {\it a priori} error bound~\eqref{eqn:a-priori-error-bound} could yield a rate of convergence of the following form: 
\begin{eqnarray}
    \| \bu - \bur \|
    = \cO(h^{p_1})
    + \cO(\Delta t^{p_2})
    + \cO\left(\left(\sum_{i=r+1}^{R} \lambda_i\right)^{p_3}\right),
    \label{eqn:rate-of-convergence}
\end{eqnarray}
where the parameters $p_1, p_2$, and $p_3$ represent the error rates of convergence for the spatial, time, and ROM discretization error components, respectively.

\item The rate of convergence~\eqref{eqn:rate-of-convergence} can yield {\it parameter scalings}, i.e., relationships between the FOM parameters and ROM parameters.
For example, one can choose the ROM dimension, $r$, in~\eqref{eqn:rate-of-convergence} so that the ROM discretization error component (i.e., the third term on the RHS of~\eqref{eqn:rate-of-convergence}) balances the spatial and time discretization error components (i.e., the first and second terms on the RHS of~\eqref{eqn:rate-of-convergence}).
This clearly illustrates the {\it practical importance} of parameter scalings:
The \apr \ guide practical decisions, e.g., choosing the ROM dimension to match a desired accuracy level.
We also note that these parameter scalings are {\it robust}, i.e., they hold over a range of problems and do not have to be recalibrated for each problem (as is commonly done in reduced order modeling).

%
%
\item We highlight that one limitation of {\it a priori} error analysis is that the error bounds rely on the exact solution, which is typically unknown in practice. 
This assumption generally does not influence the asymptotic behavior of the error bounds, as any dependence on the exact solution is often absorbed into the generic constant, $C$. However, it is important to note that the 
\apr \ cannot be explicitly used to compute accurate error estimates. In other words, the ROM lacks the certification typical of \apo, which we briefly summarize in Section~\ref{sec:a-posteriori-error-bounds}.
\end{enumerate}

\subsection{{\it A Posteriori} Error Bounds}
    \label{sec:a-posteriori-error-bounds}

The ROM \apo \ are generally of the following form:
\begin{eqnarray}
    \| \bu_h - \bur \|
    \leq s(\Delta t, h) \ d \left(\Delta t, h, \bur \right).
    \label{eqn:a-posteriori-error-bound}
\end{eqnarray}
The coefficient $s(\Delta t, h)$ usually includes the evaluation of a lower bound of an appropriate stability factor, for instance (the coercivity constant for an elliptic problem or) a \emph{inf-sup} constant for the NSE. Hence, it will normally depend on the FOM discretization parameters, but not on the ROM ones.
In contrast, the factor $d \left(\Delta t, h, \bur \right)$ shows dependence on both FOM and ROM discretization parameters. Typically, $d$ involves the evaluation of an appropriate dual norm of the residual obtained when plugging the ROM solution $\bur$ in the weak formulation of the FOM discretization, which depends on $\Delta t$ and $h$.
For more details, the interested reader is referred to, e.g., \cite{ballarin2020certified,Deparis2008,DeparisRozza2009,manzoni2014efficient}.

When comparing the \apr\ \eqref{eqn:a-priori-error-bound} to the \apo\ \eqref{eqn:a-posteriori-error-bound}, we note the following:
\begin{enumerate}
\item The quantity of interest on the left-hand side of \eqref{eqn:a-priori-error-bound} is $\| \bu - \bur \|$, while \eqref{eqn:a-posteriori-error-bound} has $\| \bu_h - \bur \|$, and hence assumes that $h$ and $\Delta t$ are fixed. In other words, the {\it a posteriori} error bound \eqref{eqn:a-posteriori-error-bound} can be used in the asymptotic regime of the ROM discretization parameter $r$ (i.e., as $r \to R$), but not in the asymptotic regime of the FOM discretization parameters $h$ and $\Delta t$.
\item In contrast to the \apr\ \eqref{eqn:a-priori-error-bound}, the \apo\ \eqref{eqn:a-posteriori-error-bound} are actually computable once the FOM discretization parameters $h$ and $\Delta t$ are fixed, 
since they do not depend on the exact solution, $\bu$.
\end{enumerate}



\subsection{A Qualitative Comparison Between {\it A Priori} and {\it A Posteriori} Error Bounds}

After 
presenting the \apr\ (Section~\ref{sec:a-priori-error-bounds}) and \apo\ (Section~\ref{sec:a-posteriori-error-bounds}), a natural question arises: 
\begin{quote}
\normalsize\it\centering
Which type of error bound should a ROM practitioner employ to make decisions
(e.g., choose discretization parameters) or ensure that given criteria are met (e.g., certify the approximations) 
in fluid flow simulations?
\end{quote}

The answer depends on the type of results one is interested in:
\begin{itemize} \itemsep5pt
	\item
        If one is interested in leveraging numerical analysis to determine {\it  parameter
        scalings} (i.e., scalings that hold across a range of problems), then \apr\ are appropriate.
        For example, as explained in Section~\ref{sec:a-priori-error-bounds}, given a FOM and a target accuracy, one could use the rate of convergence~\eqref{eqn:rate-of-convergence} to determine the appropriate ROM dimension, $r$.
        We also emphasize that these parameter scalings are robust, i.e., can be used across a range of problems.
        Thus, the expensive parameter recalibration process commonly used in reduced order modeling is avoided.        
	\item If one is interested in {\it ROM certifiability}, i.e., ensuring that the error stays within the bounds provided by the error analysis, then \apo\ are appropriate.
\end{itemize}

\section{Practical Importance of {\it A Priori} Error Bounds}
    \label{sec:practical-importance}

\subsection{Classical Methods: The ``Big O'' Notation}


We highlight that most numerical analysis for classical  methods applied to PDEs (such as finite element or finite difference methods) relies on {\it a priori} error analysis. This analysis generally provides asymptotic results resulting in error bounds which are often expressed in terms of $\cO(h^{p_1})$ 
and $\cO(\Delta t^{p_2})$. 
Furthermore, the exponent $m$ in the computational cost, $\cO(h^{m})$, is a key factor in determining the method's effectiveness.
Consequently, in the numerical analysis of classical PDE methods, the ``big O'' notation, which is an inherently asymptotic concept, plays a crucial role.


We also 
emphasize that the importance of the above notation extends beyond the numerical analysis of traditional methods for PDEs. 
For instance, in the field of numerical linear algebra~\cite{trefethen2022numerical}, computational costs (such as time and memory) are similarly described using ``big O'' notation. In this context, these costs are related to asymptotic quantities expressed as functions of the number of unknowns.

\subsection{ROM}

The central question we aim to answer in this section is
\begin{quote}
\normalsize\it\centering
    Why should a ROM practioner care about asymptotics?
\end{quote}

For FEM practitioners, it is quite natural to realize that \apr\ are a relevant numerical analysis tool, since in realistic FEM applications the number of degrees of freedom is $n = \cO(10^{6})$, or even higher. However, ROM practitioners are typically used to $r = \cO(10)$, and hence they may immediately 
question the value in exploiting \apr\ to know what happens when $r \longrightarrow \infty$.




In order to answer to the question above, and discuss our interest in the asymptotic regime, we put forth the  following two motivations:
\begin{enumerate}
    \item[(a)] The implication 
\begin{equation}
\text{\it``ROMs are low-dimensional systems of equations $\Rightarrow$ $r=\cO(10)$''} \tag{a}\label{statement1}
\end{equation}
can be 
inaccurate in general.
    \item[(b)] We consider the misconception 
    associated with the (false) implication above as one of the leading reasons for the relative underdevelopment of 
    \apr\ for ROMs.
\end{enumerate}


Below, we clarify the two points mentioned earlier, by addressing first the misconception in \eqref{statement1}.


In the authors' opinion, the primary source of this misunderstanding stems from the fact that, over the past few decades, ROMs have largely been developed within the academic community. Specifically, the numerical analysis of ROMs has been primarily conducted by the mathematical ROM community. Consequently, the test problems used to investigate the numerical analysis results were predominantly {\it academic test problems}, such as 1D or 2D Poisson equations and heat equations.

%
Even when ROMs were investigated in fluid flows, the test cases were typically simplistic academic problems. For instance, one of the most commonly used test problems is the 2D flow past a cylinder at Reynolds number $Re \leq 100$. This problem has been the workhorse for testing ROMs in fluid dynamics. The reason is straightforward: a small number of ROM basis functions can adequately capture the simple dynamics of this scenario. 
Although ROMs have been applied to fluid dynamics, they have only been successful in simulating {\it laminar} flows, which have limited practical significance.
In the authors' opinion, the next generation of ROMs in fluid dynamics should set the goal of addressing more challenging problems with real-world relevance, especially with the objective of achieving success in simulating realistic, turbulent flows.
%
%
With this objective in mind, the authors expect that $r=\cO(10)$ will not be enough to obtain accurate ROM results, hence rendering the right-hand side of \eqref{statement1} irrelevant in practice.
%

We next turn to invalidating the left-hand side of \eqref{statement1}. In simple problems, the number of degrees of freedom $n$ required by the FOM is typically moderate. For instance, a finite element discretization of a 2D flow past a cylinder at $Re = 100$ typically requires $n=\cO(10^{5})$. Assuming that a 10-dimensional ROM (i.e., with $r=\cO(10)$) can capture the large flow features, the ROM achieves a {\it four orders of magnitude} reduction in the number of degrees of freedom. However, in real-world engineering contexts, such as the nuclear engineering simulations conducted by Paul Fischer's group~\cite{fischer2022nekrs}, $n$ can be as large as $\cO(10^{8})$. Consequently, a four-order magnitude reduction would produce a ROM with a dimension of $r=\cO(10^{4})$. 
With such a value of $r$, one could not say that the ROM is a low-order model. Thus, for realistic flows, it is evident that {\it ROMs are not necessarily low-order models}.


Invalidating the relevance of both sides of \eqref{statement1} in realistic flows underscores the importance of asymptotic analysis for ROMs, just as 
in the case of classical numerical discretizations. 
We anticipate that the asymptotic analysis of ROMs will grow increasingly important in the coming years,  overcoming the limitation highlighted in (b), especially with the advent of {\it exascale computing} and 
novel grand challenge applications, where the number of degrees of freedom will increase dramatically. To address these challenges, the development of innovative {\it high-performance computing (HPC) ROMs} will be essential. For these  HPC ROMs, a numerical analysis focusing on asymptotic behavior will be crucial, much like it is for classical numerical discretizations.
Furthermore, numerical strategies that have been central in classical numerical methods will also start to become popular in reduced order modeling (see, e.g., the ROM preconditiong approaches proposed in~\cite{elman2015preconditioning,lindsay2022preconditioned}).

\section{
Contributions to {\it A Priori} Error Bounds}
    \label{sec:major-contributions}
    

In Table \ref{tab:rom_references}, we summarize 
some of the 
significant contributions to \apr\ in POD-based ROMs. 
We emphasize that this list is limited by the available space in the chapter contribution on one hand, and by the authors' personal preferences and experience on the other hand. 
We think, however, that this list still reflects some of the key  numerical analysis results addressing ROM stability, asymptotic behavior, and convergence rates. We also list contributions on the numerical analysis of ROM closures and stabilizations, which are crucial for developing ROMs applicable to realistic turbulent flows. Additionally, we highlight contributions to the numerical analysis of FOM-ROM coupling, which we believe is essential for determining {\it parameter scalings}. The references in the table are sorted choronologically, while they will be grouped by topic (last column of the table) in the subsections below.

\begin{table}
\centering
\begin{tabular}{|p{3cm}|p{1.2cm}|p{1.2cm}|p{4.5cm}|p{1.5cm}|}
\hline
\textbf{Authors} & \textbf{Year} & \textbf{Ref.} & \textbf{Description} & \textbf{Section} \\
\hhline{|=|=|=|=|=|}
Kunisch and Volkwein & 2001, 2002 & \cite{KV01, KV02} & First {\it a priori} error bounds & Section~\ref{sec:error-bounds} \\
\hline
Rathinam and Petzold & 2003 & \cite{rathinam2003new} & Alternative {\it a priori} error bounds \newline Pointwise in time error bounds & Section~\ref{sec:error-bounds} \newline Section~\ref{sec:optimal} \\
\hline
Luo et al. & 2009 & \cite{luo2008mixed} & Alternative {\it a priori} error bounds \newline Included FEM contribution \newline Pressure & Section~\ref{sec:error-bounds} \newline Section~\ref{sec:coupling} \newline Section~\ref{sec:pressure} \\
\hline
Kalashnikova and\newline Barone & 2010 & \cite{kalashnikova2010stability} & Energy-stability-preserving ROM & Section~\ref{sec:stabilization} \\
\hline
Borggaard et al. & 2011 & \cite{borggaard2011artificial} & ROM closures & Section~\ref{sec:closure} \\
\hline
Chapelle et al. & 2012 & \cite{chapelle2012galerkin} & Pointwise in time error bounds & Section~\ref{sec:optimal} \\
\hline
Iliescu and Wang & 2013 & \cite{iliescu2013variational} & Included FEM contribution & Section~\ref{sec:coupling} \\
\hline
Roop & 2013 & \cite{roop2013proper} & Closure & Section~\ref{sec:closure} \\
\hline
Iliescu and Wang & 2014 & \cite{iliescu2014are} & Pointwise in time error bounds & Section~\ref{sec:optimal} \\
\hline
Iliescu and Wang & 2014 & \cite{iliescu2014variational} & Pointwise in time error bounds \newline Included FEM contribution & Section~\ref{sec:optimal} \newline Section~\ref{sec:coupling} \\
\hline
Singler & 2014 & \cite{singler2014new} & Improved {\it a priori} error bounds & Section~\ref{sec:error-bounds} \\
\hline
Giere et al. & 2015 & \cite{giere2015supg} & FEM-ROM parameter scalings & Section~\ref{sec:coupling} \\
\hline
Kostova et al. & 2015 & \cite{kostova2015error} & Pointwise in time error bounds & Section~\ref{sec:optimal} \\
\hline
Azaiez et al & 2017 & \cite{azaiez2017streamline} & Streamline derivative ROM & Section~\ref{sec:stabilization} \\
\hline
Eroglu et al. & 2017 & \cite{eroglu2017modular} & Closure & Section~\ref{sec:closure} \\
\hline
Xie et al. & 2018 & \cite{xie2018numerical} & Stabilization with ROM filters\newline Leray ROM & Section~\ref{sec:stabilization}\newline Section~\ref{sec:closure} \\
\hline
Zerfas et al. & 2019 & \cite{zerfas2019continuous} & ROM data assimilation & Section \ref{sec:control} \\
\hline
DeCaria et al. & 2020 & \cite{decaria2020artificial} & Pressure & Section~\ref{sec:pressure} \\
\hline
Rubino & 2020 & \cite{rubino2020numerical} & Pressure & Section~\ref{sec:pressure} \\
\hline
Koc et al. & 2021 & \cite{koc2021optimal} & First optimal pointwise error bounds & Section~\ref{sec:optimal} \\
\hline
Novo and Rubino & 2021 & \cite{novo2021error} & Pressure & Section~\ref{sec:pressure} \\
\hline
Chac{\'o}n Rebollo et al. & 2022 & \cite{rebollo2022error} & Pressure & Section~\ref{sec:pressure} \\
\hline
Garc{\'\i}a-Archilla et al. & 2022 & \cite{garcia2022error} & ROM data assimilation & Section \ref{sec:control} \\
\hline
Koc et al. & 2022 & \cite{koc2022verifiability} & Data-driven ROM closures & Section~\ref{sec:closure} \\
\hline
de Frutos et al. & 2023 & \cite{de2023optimal} & ROM data assimilation & Section \ref{sec:control} \\
\hline
Eskew and Singler & 2023 & \cite{eskew2023new} & Pointwise in time error bounds & Section~\ref{sec:optimal} \\
\hline
Garc{\'\i}a-Archilla et al. & 2023 & \cite{garcia2023second} & Pointwise in time error bounds & Section~\ref{sec:optimal} \\
\hline
Garc{\'\i}a-Archilla et al. & 2023 & \cite{garcia2023pod} & Pressure & Section~\ref{sec:pressure} \\
\hline
Strazzullo et al. & 2023 & \cite{strazzullo2023new} & ROM control & Section \ref{sec:control} \\
\hline
Garc{\'\i}a-Archilla and\newline Novo & 2024 & \cite{garcia2024pointwise} & Pointwise in time error bounds & Section~\ref{sec:optimal} \\
\hline
Aza{\"\i}ez et al. & 2024 & \cite{azaiez2024least} & Pressure & Section~\ref{sec:pressure} \\
\hline
Moore et al. & 2024 & \cite{moore2024numerical} & Stabilization with ROM filters\newline Approximate deconvolution ROM & Section~\ref{sec:stabilization}\newline Section~\ref{sec:closure} \\
\hline
Reyes et al. & 2024 & \cite{reyes2024apriori} & 
Spectral element ROM parameter 
\newline 
scalings \newline  
Stabilization with ROM filters\newline Time-relaxation ROM & Section~\ref{sec:coupling} \newline \newline
Section~\ref{sec:stabilization}\newline Section~\ref{sec:closure} \\
\hline
\end{tabular}
\caption{A non-exhaustive list of contributions to {\it a priori} error bounds in POD-based ROMs.}
\label{tab:rom_references}
\end{table}

\subsection{Error Bounds}
    \label{sec:error-bounds}
%

The numerical analysis of POD-based ROMs for general problems (not limited to fluid flow) began with the seminal work of Kunisch and Volkwein, published in 2001~\cite{KV01}. In this seminal paper, the authors established the first \apr\ for ROMs, addressing both linear (e.g., heat equation) and nonlinear (i.e., Burgers equation) parabolic problems. They introduced the concept of difference quotients, 
proved inverse estimates for POD (see \cite[Lemma 2]{KV01}), and derived error bounds for the Ritz projection (see \cite[Lemma 3]{KV01}). Building on their groundbreaking work in~\cite{KV01}, the same authors addressed the more complex Navier-Stokes equations one year later~\cite{KV02}. Relevant work was also performed in~\cite{luo2008mixed,rathinam2003new}.

It took over a decade for the next significant advancement in the numerical analysis of POD-based ROMs to appear. In 2014, Singler published the paper~\cite{singler2014new} that improved the error bounds established in~\cite{KV01}. Specifically, in~\cite{singler2014new}, Singler refined the bound for the $L^2$ norm of the gradient of the POD projection error originally proven in~\cite{KV01}.

\subsection{Optimal Pointwise Rates of Convergence}
    \label{sec:optimal}

The error bounds proved in Kunisch and Volkwein's pioneering work~\cite{KV01} were for the {\it average} error.
A natural question is whether {\it pointwise} in time error bounds can also be proved.
We note that pointwise in time error bounds are standard in the numerical analysis of classical numerical methods (e.g., the finite element method~\cite{thomee2006galerkin}).
Thus, pointwise in time error bounds for POD-ROMs are also desirable.

Pointwise in time error bounds for POD-based ROMs were proved in 2014 by Iliescu and Wang~\cite{iliescu2014variational} (see~\cite{kostova2015error,rathinam2003new} for alternative bounds). 
We emphasize, however, that the pointwise in time error bounds in~\cite{iliescu2014variational} relied on~\cite[Assumption 3.1]{iliescu2014variational} (see also~\cite[Remark 3.2]{iliescu2014variational}), which was subsequently used in proving pointwise in time error bounds for parabolic problems.

A natural question is whether~\cite[Assumption 3.1]{iliescu2014variational} is essential in proving pointwise error bounds for POD-ROMs.
It turns out that this question is closely related to the question of {\it optimality} of the error bound.
This was recognized from the very beginning of the numerical analysis for POD-based ROMs.
In~\cite[Remark 1]{KV01}, Kunisch and Volkwein introduced the {\it difference quotients (DQs)} for time discretization optimality (see~\cite{chapelle2012galerkin} for an alternative strategy).
The ROM discretization optimality was first investigated in~\cite{iliescu2014are}.
We note, however, that although~\cite{iliescu2014are,KV01} suggest that DQs are needed for optimal error bounds, this statement was not rigorously proved.

The {\it first rigorous proof} of optimal pointwise in time error bounds was published in~\cite{koc2021optimal}
by Koc et al. in 2021. Specifically, the authors proved that, without using DQs, both the POD projection error and the ROM error are suboptimal not only with respect to the ROM discretization (as shown in~\cite{iliescu2014are}), but also in relation to the time discretization (as suggested in~\cite{KV01}). 
In particular, the authors constructed two analytical examples and proved that the pointwise POD projection error assumption (i.e.,~\cite[Assumption 3.1]{iliescu2014variational}) can fail, leading to a degradation by a factor of $\frac{1}{\Delta t}$ in the error bound scaling with respect to the ROM discretization error.
Furthermore, the authors proved that~\cite[Assumption 3.1]{iliescu2014variational} is automatically satisfied when the DQs are employed. 
This new finding allowed them to demonstrate that, with the DQs, both the POD projection error and the ROM error achieve optimal convergence rates.
Finally, the authors revisited the definition of ROM discretization error optimality by introducing a stronger notion of optimality. 
This novel definition allowed them to show that all the pointwise in time error bounds in the DQ case are optimal in at least one sense.

The new results in~\cite{koc2021optimal} have spurred several new research developments.
In~\cite{eskew2023new,garcia2023second}, it was shown that, to prove optimal pointwise error bounds, it is sufficient to use only the DQs and either the snapshot at the initial time or the mean value of the snapshots.
Further improvements were recently proved in~\cite{garcia2024pointwise}.

\subsection{FOM-ROM Coupling and Parameter Scalings}
\label{sec:coupling}

We highlight that 
the numerical analysis of ROMs for fluids often ignores the FOM discretization used to generate the snapshots. As a result, these ROM error bounds 
{\it do not depend on the FOM parameters}.


The typical rationale for disregarding the FOM numerical discretization can be illustrated for the 
finite element method. 
In that case, one can argue that the discretization 
essentially functions as a weighting matrix in the snapshot generation process~\cite[Chapter 1, Section 3]{volkwein2013proper}. While this reasoning is technically accurate for simple PDEs, such as the Poisson or heat equations discussed in~\cite{volkwein2013proper}, the situation is more complicated for the Navier-Stokes equations, which involve a nonlinear term and two different types of variables (velocity and pressure), typically requiring different finite elements (i.e., mixed methods).


The first numerical analysis of POD-based ROMs which explicitly considers the FOM discretization was conducted by Luo et al. in 2009~\cite{luo2008mixed}.
A different strategy was used by Iliescu and Wang in 2013, 2014~\cite{iliescu2013variational,iliescu2014variational}, where the authors established error bounds for a stabilized ROM for the convection-diffusion-reaction equation and Navier-Stokes equations. 
In these papers, the error bounds explicitly depend on the finite element discretization used at the FOM level to generate the snapshots.


We stress that deriving error bounds that depend on both ROM and FOM parameters 
allows us to establish {\it rigorous parameter scalings} for the ROMs, similar to what is typically done for classical numerical discretizations (see the monograph~\cite{roos2008robust}). In contrast, error bounds that disregard the underlying FOM discretization~\cite{volkwein2013proper} prevent the determination of these critical ROM parameter scalings. 
The first rigorous parameter scalings were proven in~\cite{giere2015supg} by Giere et al. (2015), who used error bounds dependent on both FEM and ROM parameters to establish parameter scalings for the streamline-upwind Petrov-Galerkin (SUPG) ROM stabilization strategy.
In~\cite{reyes2024apriori}, Reyes et al. (2024) leveraged \apr\ to prove robust parameter scalings for the time-relaxation ROM, which is a filter-based stabilization (regularization).
Specifically, the authors derived scalings between the time-relaxation parameter and the filter radius and ROM dimension.

\subsection{Pressure}
\label{sec:pressure}

While the crucial role of pressure in ROMs was acknowledged early on (see, for instance, Noack et al. in 2005~\cite{noack2005need}, as well as~\cite{bergmann2009enablers,caiazzo2014numerical,weller2010numerical}), most of the 
works on \apr\ neglect to approximate the pressure at the ROM level. 
Luo~et al.~\cite{luo2008mixed} have presented \apr\ for the ROM pressure of the time-dependent Navier-Stokes equations discretized by the finite element method. 
In 2020, two significant developments on the numerical analysis of ROM stabilizations for pressure were published almost simultaneously: DeCaria et al.~\cite{decaria2020artificial} established error bounds for the artificial compression ROM, while Rubino~\cite{rubino2020numerical} provided error bounds for the local projection stabilization. Since then, additional contributions have advanced the numerical analysis of ROM pressure.
Novo and Rubino~\cite{novo2021error} and Garc{\'\i}a-Archilla et al.~\cite{garcia2023pod} used grad-div stabilization at both the FOM and the ROM levels, and proved error bounds for the pressure that are independent of inverse powers of the viscosity.
Chac{\'o}n Rebollo et al.~\cite{rebollo2022error} proved \apr\ for the pressure in a residual-based stabilization-motivated POD-ROM.
Azaiez et al.~\cite{azaiez2024least} approximated the pressure by solving a least-squares problem for the residual of the reduced velocity with respect to a dual norm, and proved \apr\ for the pressure.

\subsection{Stabilization}
\label{sec:stabilization}
A summary of existing numerical analyses for residual-based stabilizations of the convection-diffusion-reaction equation can be found in~\cite[Section 6.3]{parish2023residual}.
The \apr\ for POD-based ROM stabilizations of convection-dominated fluid flows are relatively scarce.
For example, Kalashnikova and
Barone~\cite{kalashnikova2010stability} derived \apr\ for an
energy-stability-preserving ROM formulation introduced in
\cite{barone2009stable} for linearized compressible flow.  These error bounds employed a carefully constructed stable penalty-like implementation of the relevant boundary conditions in the ROM.
Azaiez et al.~\cite{azaiez2017streamline} proved \apr\ for a streamline derivative POD-ROM stabilization for convection-dominated flows.
Xie et al.~\cite{xie2018numerical} analyzed the Leray ROM, which is a filter-based ROM regularization.
In particular, the authors proved the first error bounds for the ROM differential filter, which they leveraged to prove \apr\ for the Leray ROM.
Recently, Moore et al.~\cite{moore2024numerical} proved the first error bounds for approximate deconvolution and \apr\ for the approximate deconvolution Leray ROM.
Finally, Reyes et al.~\cite{reyes2024apriori} proved \apr\ for the time-relaxation ROM, which is a filter-based ROM stabilization for convection-dominated flows.

\subsection{Closure}
\label{sec:closure}

The first investigations on numerical analysis for ROM closures were conducted by Borggaard et al. in 2011~\cite{borggaard2011artificial}. In their paper, the authors examined the 1D Smagorinsky ROM closure model, demonstrating its unconditional stability and convergence with respect to time discretization and POD truncation~\cite[Theorem 4.1]{borggaard2011artificial}. Subsequent numerical analyses of various ROM closures were covered in a series of papers recently reviewed in~\cite{ahmed2021closures}.
The numerical analysis of variational multiscale ROMs was performed by Iliescu and Wang in~\cite{iliescu2014variational}, where \apr\  were rigorously proven. 
Further developments were presented by Roop~\cite{roop2013proper} and Eroglu et al.~\cite{eroglu2017modular}.
We also refer to Rubino~\cite{rubino2020numerical} for numerical analysis for ROM stabilizations of convection-dominated flows.

In 2018, Xie et al. 
carried out the numerical analysis of ROM filters~\cite{xie2018numerical}, proving error bounds for the ROM differential filter~\cite[Lemma 4.3]{xie2018numerical}. Furthermore, in 2022 Koc et al. performed the 
numerical analysis of data-driven ROM closures~\cite{koc2022verifiability}, proving the verifiability of a new data-driven ROM closure and showing that the ROM approximation improves upon improving the ROM closure approximation~\cite[Theorem 2]{koc2022verifiability}. 
Recently, 
Moore et al. conducted the first numerical analysis of ROM approximate deconvolution~\cite{moore2024numerical}, and Reyes et al.~\cite{reyes2024apriori}  proved \apr\ for the time-relaxation ROM.

\subsection{Data Assimilation and Control}
\label{sec:control}

Although ROMs have been used in data assimilation \cite{karcher2018reduced,maday2015parameterized,maday2013generalized} and control \cite{alla2018posteriori,hinze2005proper,negri2015reduced,negri2013reduced,strazzullo2018model,strazzullo2020pod}, the \apr\ are relatively scarce.

Zerfas et al.~\cite{zerfas2019continuous} proposed a nudging data assimilation algorithm to improve the ROM long-time accuracy.
The authors proved that, with a properly chosen nudging parameter, the proposed nudging algorithm converges exponentially fast in time to the true solution. 
Furthermore, the authors proposed a strategy for nudging adaptively in time, by adjusting dissipation arising from the nudging term to better match true solution energy.
A nudging algorithm for data assimilation was also considered by Garc{\'\i}a-Archilla et al.~\cite{garcia2022error}, who used grad-div stabilization at the ROM level and proved \apr\ with constants independent on inverse powers of the viscosity.

In~\cite{de2023optimal}, de Frutos et al. considered the reduced order modeling of infinite horizon problems via the dynamic programming approach, and carried out the error analysis of the method.
Strazzullo et al.~\cite{strazzullo2023new} proposed a new feedback control strategy for high Reynolds number flows, and proved \apr\ for high Reynolds numbers that were not covered by current results.

\section{Conclusions and Outlook}
    \label{sec:conclusions}

In this paper, we presented a brief survey of \apr\ for POD-based ROMs of fluid flows.
We outlined the main components of \apr\ and included a short qualitative comparison with \apo.
We also emphasized the importance of \apr\ in practical settings, e.g., when the ROM dimension is inherently large due to the complex dynamics of the underlying system and the large FOM dimension (e.g., $\cO(10^8)$).
Finally, we summarized some of the contributions to \apr\ for POD-based ROMs, including the first proofs for fundamental results, optimal pointwise rates of convergence, FOM-ROM coupling and parameter scalings, \apr\ for the pressure, stabilizations, and closures, as well as \apr\ in data assimilation and control.

We tried to cover as many developments as possible given the chapter contribution page limitations.
Despite our efforts, we are aware that important contributions have not been discussed.
We believe, however, that this brief survey could serve as a stepping stone toward a comprehensive review of the exciting research area of numerical analysis for ROMs of fluid flows.

\begin{acknowledgement}
We acknowledge the European Union's Horizon 2020 research and innovation program under the Marie Skłodowska-Curie Actions, grant agreement 872442 (ARIA).
FB thanks the PRIN 2022 PNRR project ``ROMEU: Reduced Order Models for Environmental and Urban flows'' (CUP J53D23015960001). TI acknowledges support through National Science Foundation grants DMS-2012253 and CDS\&E-MSS-1953113.
\end{acknowledgement}

\ethics{Competing Interests}{ 
The authors have no conflicts of interest to declare that are relevant to the content of this chapter.}

\bibliographystyle{plain}
\bibliography{traian}

\end{document}